\documentstyle[11pt,leqno]{article}

\newtheorem{theorem}{Theorem}[section]

\newtheorem{lemma}[theorem]{Lemma}

\newtheorem{example}[theorem]{Example}
\newtheorem{observation}[theorem]{Observation}

\newcommand{\pf}{ {\em Proof.\ \ }}

\def\eop{\hfill\qquad\rule[-1mm]{1.75mm}{1.75mm}}

\newcommand{\be}{\begin{equation}}
\newcommand{\ee}{\end{equation}}

\newcommand{\lm}{\lambda=\frac{\mu}{1+\mu}}

\newcommand{\ts}{\tau_s=\frac{\sigma_s}{1+\sigma_s}}

\newcommand{\G}{{\cal G}}
\newcommand{\RG}{{\cal R}}

\newcommand{\ls}{L_s}

\newcommand{\R}{\mbox{\bf R}}
\newcommand{\CO}{\mbox{\bf C}}

\title{Z-Pencils}

\author{
J. J. McDonald \footnotemark[1]\ \ \footnotemark[9]
\and
D. D. Olesky \footnotemark[2]\ \ \footnotemark[9]
\and
H. Schneider \footnotemark[3]
\and
M. J. Tsatsomeros \footnotemark[1]\ \ \footnotemark[9]
\and
P. van den Driessche \footnotemark[4]\ \ \footnotemark[9]
}

\begin{document}
\maketitle

\renewcommand{\thefootnote}{\fnsymbol{footnote}}

\footnotetext[1]{Dept of Mathematics and Statistics, Univ. of
Regina, Regina, Saskatchewan S4S 0A2.}

\footnotetext[2]{Dept of Computer Science,
Univ. of Victoria, Victoria, British Columbia V8W 3P6.}

\footnotetext[3]{Dept of Mathematics, Univ. of
Wisconsin, Madison, Wisconsin 53706.}

\footnotetext[4]{Dept of Mathematics and Statistics,
Univ. of Victoria, Victoria, British Columbia V8W 3P4.}

\footnotetext[9]{Research partially supported by NSERC research grant.}

\renewcommand{\thefootnote}{\arabic{footnote}}

\begin{abstract}
The matrix pencil $(A,B)=\{tB-A \ | \ t\in\CO\}$ is considered under
the assumptions that $A$ is entrywise nonnegative and $B-A$ is a nonsingular
M-matrix. As $t$ varies in $[0,1]$, the Z-matrices $tB-A$ are partitioned 
into the sets $\ls$ introduced by Fiedler and Markham. 
As no combinatorial structure of $B$ is assumed here, 
this partition generalizes some of their work where $B=I$.
Based on the union of the directed graphs of $A$ and $B$, the combinatorial
structure of nonnegative eigenvectors associated with the largest eigenvalue of 
$(A,B)$ in $[0,1)$ is considered.
\end{abstract}



\pagestyle{myheadings}
\thispagestyle{plain}

\section{Introduction}

The generalized eigenvalue problem $Ax=\lambda Bx$ for
$A=[a_{ij}], \ B=[b_{ij}]\in\R^{n,n}$, with inequality
conditions motivated by certain economics models, was studied 
by Bapat et al. \cite{BOV}. In keeping with this work, we consider the matrix
pencil $(A,B)=\{tB-A \ | \ t\in\CO\}$ under the conditions
\be \label{c1}\mbox{$A$ is entrywise nonnegative, denoted by $A\geq 0$} \ee
\be \label{c2}\mbox{$b_{ij}\leq a_{ij}$ for all $i\neq j$} \ee
\be \label{c3}\mbox{there exists a positive vector $u$ such that $(B-A)u$
                    is positive}. \ee
Note that in \cite{BOV} $A$ is also assumed to be irreducible, but that
is not imposed here. When $Ax=\lambda Bx$ for some nonzero $x$,
the scalar $\lambda$ is an {\em eigenvalue} and $x$ is the
corresponding {\em eigenvector} of $(A,B)$. The {\em eigenspace}
of $(A,B)$ associated with an eigenvalue $\lambda$ is the nullspace
of $\lambda B-A$.

A matrix $X\in\R^{n,n}$ is a {\em Z-matrix} if $X=qI-P$, where $P\geq 0$ and
$q\in\R$. If, in addition, $q\geq\rho(P)$, where $\rho(P)$ is the spectral 
radius of $P$, then $X$ is an {\em M-matrix}, and is singular if and only if
$q=\rho(P)$. It follows from (\ref{c1}) and (\ref{c2}) that when $t\in [0,1]$,
$tB-A$ is a Z-matrix. Henceforth the term {\em Z-pencil} $(A,B)$ refers to
the circumstance that $tB-A$ is a Z-matrix for all $t\in [0,1]$.

Let $\langle n\rangle=\{1,2,\ldots,n\}$. If $J\subseteq\langle n\rangle$,
then $X_J$ denotes the principal submatrix of $X$ in rows and columns of $J$.
As in \cite{FM}, given a nonnegative $P\in\R^{n,n}$
and an $s\in\langle n\rangle$, define
\[ \rho_s(P)=\max_{|J|=s}\{\rho(P_J)\} \]
and set $\rho_{n+1}(P)=\infty$.
Let $\ls$ denote the set of Z-matrices in $\R^{n,n}$ of the form $qI-P$,
where $\rho_s(P)\leq q<\rho_{s+1}(P)$ for $s\in\langle n\rangle$,
and $-\infty<q<\rho_1(P)$ when $s=0$. This gives a partition of all
Z-matrices of order $n$. Note that $qI-P\in L_0$ if and only if
$q<p_{ii}$ for some $i$. Also, $\rho_n(P)=\rho(P)$, and $L_n$ is
the set of all (singular and nonsingular) M-matrices.

We consider the Z-pencil $(A,B)$ subject to conditions
(\ref{c1})-(\ref{c3}) and partition its matrices into the sets $\ls$. 
Viewed as a partition of the Z-matrices $tB-A$ for $t\in [0,1]$,
our result provides a generalization of some of 
the work in \cite{FM} (where $B=I$).
Indeed, since no combinatorial structure of $B$ is assumed, 
our Z-pencil partition is a consequence of a more complicated 
connection between the Perron-Frobenius theory for $A$ and the spectra 
of $tB-A$ and its submatrices. 

Conditions (\ref{c2}) and (\ref{c3}) imply that $B-A$ is a 
nonsingular M-matrix and thus its inverse is entrywise nonnegative
(see \cite[N$_{38}$, p. 137]{BP}). This, together with (\ref{c1}),
gives $(B-A)^{-1}A\geq 0$. Perron-Frobenius theory is used in \cite{BOV} to
identify an eigenvalue $\rho(A,B)$ of the pencil $(A,B)$, defined as
\[ 
\rho(A,B)=\frac{\rho\left((B-A)^{-1}A\right)}{1+\rho\left((B-A)^{-1}A\right)}. 
\]
Our partition involves $\rho(A,B)$ and the eigenvalues of the
subpencils $(A_J,B_J)$. Our Z-pencil partition result, Theorem \ref{zthm}, 
is followed by examples where as $t$ varies in $[0,1]$, 
$tB-A$ ranges through some or all of the sets $\ls$ for $0\leq s\leq n$.
In Section 3 we turn to a consideration of the combinatorial structure of
nonnegative eigenvectors associated with $\rho(A,B)$. This involves some 
digraph terminology, which we introduce at the beginning of that section.

In \cite{FM}, \cite{Sm} and \cite{Na}, interesting results on the spectra of 
matrices in $\ls$, and a classification in terms of the inverse of a Z-matrix,
are established. These results are of course applicable to the matrices of a
Z-pencil, however, as they do not directly depend on the form $tB-A$ of 
the Z-matrix, we do not consider them here.

\section{Partition of Z-pencils}

We begin with two observations and a lemma used to prove our result on the 
Z-pencil partition.
\begin{observation}\label{obs1}  
Let $(A,B)$ be a pencil with $B-A$ nonsingular.  
Given a real $\mu\neq -1$, let $\lm$. Then
the following hold:
\begin{description}
\item[(i)] $\lambda\neq 1$ is an eigenvalue of $(A,B)$ if
and only if $\mu\neq -1$ is an eigenvalue of $(B-A)^{-1}A$. 
\item[(ii)] 
$\lambda$ is a strictly increasing function of $\mu\neq -1$.
\item[(iii)] 
$\lambda \in [0,1)$ if and only if $\mu\geq 0$.
\end{description}
\end{observation}

\pf
If $\mu$ is an eigenvalue of $(B-A)^{-1}A$, then
there exists nonzero $x\in\R^n$ such that $(B-A)^{-1}Ax=\mu x$.
It follows that $Ax=\mu(B-A)x$ and if $\mu\neq -1$, then
$Ax=\frac{\mu}{1+\mu}Bx=\lambda Bx$. Notice that $\lambda$ cannot be $1$ for
any choice of $\mu$. The reverse argument shows that the
converse is also true. The last statement of (i) is obvious.
Statements (ii)  and (iii) follow easily from the definition of $\lambda$.
\eop

Note that $\lambda=1$ is an eigenvalue of $(A,B)$ if and only if 
$B-A$ is singular.

\begin{observation}\label{obs2} 
Let $(A,B)$ be a pencil satisfying (\ref{c2}), (\ref{c3}). 
Then the following hold:
\begin{description}
\item[(i)] For any nonempty $J\subseteq \langle n \rangle$,
           $B_J-A_J$ is a nonsingular M-matrix.
\item[(ii)] If in addition (\ref{c1}) holds, the largest real 
            eigenvalue of $(A,B)$ in $[0,1)$ is $\rho(A,B)$.
\end{description}
\end{observation}

\pf
(i) This follows since (\ref{c2}) and (\ref{c3}) imply that $B-A$
is a nonsingular M-matrix (see \cite[I$_{27}$, p. 136]{BP})
and since every principal submatrix of a nonsingular M-matrix
is also a nonsingular M-matrix (see \cite[p. 138]{BP}).

(ii) This follows from Observation \ref{obs1}, 
since $\mu=\rho((B-A)^{-1}A)$ is the maximal
positive eigenvalue of $(B-A)^{-1}A$.
\eop

\begin{lemma} \label{zlem} 
Let $(A,B)$ be a pencil satisfying (\ref{c1})-(\ref{c3}).
Let $\mu=\rho\left((B-A)^{-1}A\right)$ and $\rho(A,B)=\frac{\mu}{1+\mu}$. 
Then the following hold:
\begin{description}
\item[(i)] For all $t\in \left(\rho(A,B), 1\right]$, 
           $tB-A$ is a nonsingular M-matrix.
\item[(ii)] The matrix $\rho(A,B) B-A$ is a singular M-matrix.
\item[(iii)]  For all $t\in \left(0,\rho(A,B)\right)$, 
              $tB-A$ is not an M-matrix. 
\item[(iv)] For $t=0$, either $tB-A$ is a singular M-matrix or is not an
            M-matrix.
\end{description}
\end{lemma}

\pf
Recall that (\ref{c1}) and (\ref{c2}) imply that $tB-A$ is
a Z-matrix for all $0<t\leq 1$. As noted in Observation \ref{obs2} (i),
$B-A$ is a nonsingular M-matrix and thus its eigenvalues 
have positive real parts \cite[G$_{20}$, p. 135]{BP}, and the
eigenvalue with minimal real part is real (\cite[Exercise 5.4, p. 159]{BP}.
Since the eigenvalues are continuous functions of the entries
of a matrix, as $t$ decreases from $t=1$, $tB-A$ is a nonsingular M-matrix
for all $t$ until a value of $t$ is encountered for which $tB-A$ is
singular. Results (i) and (ii) now follow by Observation \ref{obs2} (ii).

To prove (iii), consider $t\in \left(0,\rho(A,B)\right)$.
Since $(B-A)^{-1}A\geq 0$, there exists an eigenvector
$x\geq 0$ such that $(B-A)^{-1}Ax=\mu x$. Then
$Ax=\rho(A,B) Bx$ and $(tB-A)x=\left(t-\rho(A,B)\right)Bx\leq 0$ since
$Bx=\frac{1}{\rho(A,B)}Ax\geq 0$.  By \cite[A$_5$,
p. 134]{BP},  $tB-A$ is not a nonsingular M-matrix. 
To complete the proof (by contradiction), suppose
$\alpha B-A$ is a singular M-matrix for some $\alpha\in
(0,\rho(A,B))$. Since there are finitely many values of $t$ for which 
$tB-A$ is singular, we can choose $\beta\in (\alpha,\rho(A,B))$ such that 
$\beta B-A$ is nonsingular. Let $\epsilon=\frac{\beta-\alpha}{\alpha}$. 
Then $(1+\epsilon)(\alpha B-A)$ is a singular M-matrix and  
\[
(1+\epsilon)(\alpha B-A)+\gamma I=\beta B-A-\epsilon A+\gamma I 
\leq \beta B-A+\gamma I  \] 
since $A\geq 0$ by (\ref{c1}). By \cite[C$_9$, p. 150]{BP}, 
$\beta B-A-\epsilon A +\gamma I$ is a nonsingular
M-matrix for all $\gamma>0,$ and hence $\beta B-A +\gamma I$
is a nonsingular M-matrix
for all $\gamma>0$ by \cite[2.5.4, p. 117]{HJ}. This implies that
$\beta B-A$ is also a (nonsingular) M-matrix 
(\cite[C$_9$, p. 150]{BP}), contradicting the above. Thus we
can also conclude that $\alpha B-A$ cannot be a singular M-matrix
for any choice of $\alpha\in (0,\rho(A,B))$, establishing (iii).
For (iv), $-A$ is a singular M-matrix if and only if it is,
up to a permutation similarity, strictly triangular. 
Otherwise, $-A$ is not an M-matrix.
\eop

\begin{theorem} \label{zthm} 
Let $(A,B)$ be a pencil satisfying (\ref{c1})-(\ref{c3}).
For $s=1,2,\ldots,n$ let 
\[
\sigma_s=\max_{|J|=s}\{\rho\left((B_J-A_J)^{-1}A_J\right)\}, 
\ \ \ \ts,  
\]
and \ $\tau_{0}=0$. Then for $s=0,1,\ldots, n-1$ and 
$\tau_s\leq t<\tau_{s+1}$, the matrix $tB-A\in \ls$.
For $s=n$ and $\tau_n\leq t\leq 1$, the matrix $tB-A\in L_n$.
\end{theorem}

\pf
Fiedler and Markham \cite[Theorem 1.3]{FM} show that
for $1\leq s\leq n-1$,
$X\in \ls$ if and only if all principal submatrices of $X$ of
order $s$ are M-matrices, and there exists a
principal submatrix of order $s+1$ that is not an M-matrix.
Consider any nonempty $J\subseteq \langle n \rangle$ and $t\in [0,1]$.
Conditions (\ref{c1}) and (\ref{c2}) imply that $tB_J-A_J$ is a Z-matrix.
By Observation \ref{obs2} (i), $B_J-A_J$ is a nonsingular M-matrix.
Let $\mu_J=\rho\left((B_J-A_J)^{-1}A_J\right)$. Then by
Observation \ref{obs2} (ii), $\tau_J=\frac{\mu_J}{1+\mu_J}$ is
the largest eigenvalue in $[0,1)$ of the pencil
$(A_J,B_J)$. Combining this with Observation \ref{obs2} (i) and Lemma 
\ref{zlem}, the matrix $tB_J-A_J$ is an M-matrix for all $\tau_J\leq
t\leq 1,$ and $tB_J-A_J$ is not an M-matrix for all
$0< t< \tau_J$. If $1\leq s\leq n-1$ and $|J|=s$, then $tB_J-A_J$ 
is an M-matrix for all $\tau_s \leq t\leq 1$. Suppose $\tau_s<\tau_{s+1}$.
Then there exists $K\subseteq\langle n\rangle$ such
that $|K|=s+1$  and $tB_K-A_K$ is not an M-matrix for $0<t<\tau_{s+1}$.
Thus by \cite[Theorem 1.3]{FM} $tB-A\in\ls$ for all $\tau_s\leq t<\tau_{s+1}$.
When $s=n$, since $B-A$ is a nonsingular M-matrix, 
$tB-A\in L_n$ for all $t$ such that $\rho(A,B)=\tau_n\leq t\leq 1$ 
by Lemma \ref{zlem} (i).
For the case $s=0$, if $0<t<\tau_1$, then $tB-A$ has a negative
diagonal entry and thus $tB-A\in L_0$. For $t=0$,
$tB-A=-A$. If $a_{ii}\neq 0$ for some $i\in\langle n \rangle,$
then $-A\in L_0$; if $a_{ii}=0$ for all $i\in\langle n \rangle$,
then $\tau_1=\tau_0=0$, namely, $-A\in L_s$ for some $s\geq 1$.
\eop

We continue with illustrative examples.
\begin{example}
\label{example1}
{\em Consider 
\[ A=\pmatrix{1&2\cr 
              1&0}\ \ \mbox{and} \ \ \
B=\pmatrix{2&2\cr 
           1&1}, \]
for which $\tau_2=2/3$ and $\tau_1=1/2$.
It follows that 
\[
tB-A\ \in \ \left\{ \begin{array}{lllll}
       L_0 & \mbox{if $0\leq t<1/2$}\\ \\
       L_1 & \mbox{if $1/2\leq t< 2/3$}\\ \\
       L_2 & \mbox{if $2/3\leq t\leq 1$.}
                   \end{array}\right.
\]
That is, as $t$ increases from $0$ to $1$, $tB-A$ belongs to all the 
possible Z-matrix classes $L_s$.
}
\end{example}

\begin{example}
\label{example2} 
{\em Consider the matrices in \cite[Example 5.3]{BOV}, that is,
\[ A=\pmatrix{1&0&0&0\cr 
              1&1&0&0\cr
              0&0&1&0\cr 
              0&1&0&1}\ \ \mbox{and} \ \ \
B=\pmatrix{4&0&-2&0\cr 
           0&3&0&-1\cr
          -2&0&4&0\cr 
           0&-2&0&4}. \]
Referring to Theorem \ref{zthm}, $\tau_4=\rho(A,B)=\frac{4+\sqrt{6}}{10}=\tau_3=\tau_2$ and $\tau_1=1/3$.
It follows that
\[
tB-A\ \in \ \left\{ \begin{array}{lllll}
           L_0 & \mbox{if $0\leq t<1/3$}\\ \\
           L_1 & \mbox{if $1/3\leq t< \frac{4+\sqrt{6}}{10}$}\\ \\
           L_4 & \mbox{if $\frac{4+\sqrt{6}}{10}\leq t\leq 1$.}
                   \end{array}\right.
\]
Notice that for $t\in [0,1]$, $tB-A$ ranges through 
only $L_0$, $L_1$ and $L_4$.
}
\end{example}

\begin{example}
\label{example3}
{\em Now let 
\[ A=\pmatrix{0&1\cr 
              0&0}\ \ \mbox{and} \ \ \
B=\pmatrix{1&1\cr  
           0&1}. \]
In contrast to the above two examples, $tB-A\in L_2$ for all $t\in [0,1]$.
Note that, in general, $tB-A\in L_n$ for all $t\in[0,1]$ if and
only if $\rho(A,B)=0$.}
\end{example}

\section{Combinatorial Structure of the Eigenspace of $\rho(A,B)$}

Let $\Gamma=(V,E)$ be a {\it digraph}, where $V$ is a finite
vertex set and $E\subseteq V\times V$ is the edge set.  
If $\Gamma'=(V,E')$, then $\Gamma\cup\Gamma'=(V,E\cup E')$.
Also write $\Gamma'\subseteq\Gamma$ when $E'\subseteq E$.
For $j\neq k$,
a {\em path of length} $m\geq 1$ from $j$ to $k$ in $\Gamma$ is a
sequence of vertices $j=r_1,r_2,\ldots,r_{m+1}=k$ such that
$(r_s,r_{s+1})\in E$ for $s=1,\ldots,m$. As in \cite[Ch. 2]{BP},
if $j=k$ or if there is a path from vertex $j$ to vertex $k$ in $\Gamma$, 
then $j$ has {\em access to} $k$ (or $k$ is {\em accessed from} $j$).  
If $j$ has access to $k$ and $k$ has access to $j$, then $j$ and $k$ 
{\em communicate}. The communication
relation is an equivalence relation, hence $V$ can be partitioned into
equivalence classes, which are referred to as the {\em classes} of $\Gamma$.

The {\em digraph} of $X=[x_{ij}]\in\R^{n,n}$, denoted by $\G(X)=(V,E)$,
consists of the vertex set $V=\langle n\rangle$ and the
set of directed edges $E=\{(j,k) \ | \ x_{jk}\ne 0\}$.
If $j$ has access to $k$ for all distinct $j,k\in V$, then
$X$ is {\em irreducible} (otherwise, {\em reducible}). 
It is well known  that the rows and 
columns of $X$ can be simultaneously reordered so that
$X$ is in block lower triangular {\em Frobenius normal form}, with 
each diagonal block irreducible. The irreducible blocks
in the Frobenius normal form of $X$ correspond 
to the classes of $\G(X)$. 

In terminology similar to that of \cite{S},
given a digraph $\Gamma$, the {\em reduced graph} of $\Gamma$,
$\RG(\Gamma)=(V',E')$, is the digraph derived from $\Gamma$ by taking
\[ V'=\{J \ |\ J \ \mbox{is a class of\ $\Gamma$}\} \] 
and
\[ E'=\{(J,K) \ | \ \mbox{there exist $j\in J$ and $k\in K$ such that
                       $j$ has access to $k$ in $\Gamma$}\}.  \]
When $\Gamma=\G(X)$ for some $X\in\R^{n,n}$, we denote $\RG(\Gamma)$ 
by $\RG(X)$.

Suppose now that $X=qI-P$ is a singular M-matrix, where $P\geq 0$ 
and $q=\rho(P)$. If an irreducible block $X_J$ in the Frobenius normal 
form of $X$ is singular, then $\rho(P_J)=q$ and we refer to the 
corresponding class $J$ as a {\em singular class} (otherwise, a 
{\em nonsingular class}). A singular class $J$ of $\G(X)$ is called 
{\em distinguished} if when $J$ is accessed from a class $K\neq J$ 
in $\RG(X)$, then $\rho(P_K)<\rho(P_J)$. That is, a singular class $J$ 
of $\G(X)$ is distinguished if and only if $J$ is accessed only from itself 
and nonsingular classes in $\RG(X)$.

We paraphrase now Theorem 3.1 of \cite{S} as follows.
\begin{theorem}
\label{pb}
Let $X\in\R^{n,n}$ be an M-matrix and let $J_1,\ldots,J_p$
denote the distinguished singular classes of $\G(X)$.
Then there exist unique (up to scalar multiples) nonnegative
vectors $x^1,\ldots,x^p$ in the nullspace of $X$ such that
\[ x^i_j\ \left\{
\begin{array}{l}
\mbox{$=0$ if $j$ does not have access to a vertex in $J_i$ in $\G(X)$} \\
\mbox{$>0$ if $j$ has access to a vertex in $J_i$ in $\G(X)$}
\end{array}\right. \]
for all $i=1,2,\ldots p$ and $j=1,2,\ldots, n$. Moreover,
every nonnegative vector in the nullspace of $X$ is a linear
combination with nonnegative coefficients of $x^1,\ldots,x^p$.
\end{theorem}
We apply the above theorem to a Z-pencil, using the following
lemma.
\begin{lemma} \label{graph}
Let $(A,B)$ be a pencil satisfying (\ref{c1}) and (\ref{c2}). Then
the classes of $\G(tB-A)$ coincide with the classes of 
$\G(A)\cup \G(B)$ for all $t\in (0,1)$.
\end{lemma}

\pf
Clearly $\G(tB-A)\subseteq\G(A)\cup\G(B)$ for all scalars $t$. 
For any $i\neq j$, if either $b_{ij}\neq 0$ or $a_{ij}\neq 0$, 
and if $t\in (0,1)$, conditions (\ref{c1}) and (\ref{c2}) imply that
$tb_{ij}<a_{ij}$ and hence $tb_{ij}-a_{ij}\neq 0$. This means that apart 
from vertex loops, the edge sets of $\G(tB-A)$ and $\G(A)\cup \G(B)$ 
coincide for all $t\in (0,1)$.
\eop

\begin{theorem}
\label{zbasis} 
Let $(A,B)$ be a pencil satisfying (\ref{c1})-(\ref{c3}) and let
\[
\Gamma=\left\{ \begin{array}{lll}
           \G(A)\cup\G(B) & \mbox{if $\rho(A,B)\neq 0$}\\ \\
           \G(A) & \mbox{if $\rho(A,B)=0$.}
                   \end{array}\right.
\]
Let $J_1,\ldots,J_p$ denote the classes of 
$\Gamma$ such that for each $i=1,2,\ldots,p$,

\hspace*{.25in}(i) \ $(\rho(A,B)B-A)_{J_i}$ is singular, and 

\hspace*{.25in}(ii)\ if $J_i$ is accessed from a class
$K\neq J_i$ in $\RG(\Gamma)$, then $(\rho(A,B)B-A)_K$ is nonsingular.

Then there exist unique (up to scalar multiples) nonnegative
vectors $x^1,\ldots,x^p$ in the eigenspace associated with the 
eigenvalue $\rho(A,B)$ of $(A,B)$ such that
\[
x^i_j\ \left\{ 
\begin{array}{l}
\mbox{$=0$ if $j$ does not have access to a vertex in $J_i$ in $\Gamma$} \\
\mbox{$>0$ if $j$ has access to a vertex in $J_i$ in $\Gamma$}
\end{array}\right.
\]
for all $i=1,2,\ldots, p$ and $j=1,2,\ldots, n$.
Moreover, every nonnegative vector in the eigenspace associated with the 
eigenvalue $\rho(A,B)$ is a linear combination with nonnegative 
coefficients of $x^1,\ldots,x^p$.

\end{theorem}

\pf
By Lemma \ref{zlem} (ii), $\rho(A,B)B-A$ is a singular M-matrix. Thus
\[ \rho(A,B)B-A=qI-P=X, \]
where $P\geq 0$ and $q=\rho(P)$. 
When $\rho(A,B)=0$, the result follows from Theorem \ref{pb} applied
to $X=-A$. When $\rho(A,B)>0$, by Lemma \ref{graph},
$\Gamma=\G(X)$. Class $J$ of $\Gamma$ is singular
if and only if $\rho(P_J)=q$, which is equivalent
to $(\rho(A,B)B-A)_{J}$ being singular. Also a singular 
class $J$ is distinguished
if and only if for all classes $K\neq J$ that access $J$ in $\RG(X)$,
$\rho(P_K)<\rho(P_J)$, or equivalently $(\rho(A,B)B-A)_{K}$ is nonsingular.
Applying Theorem \ref{pb} gives the result. 
\eop

We conclude with a generalization of Theorem 1.7 of \cite{FM} to Z-pencils.
Note that the class $J$ in the following result is a singular
class of $\G(A)\cup \G(B)$.
\begin{theorem}
\label{zs}
Let $(A,B)$ be a pencil satisfying (\ref{c1})-(\ref{c3})
and let $t\in(0,\rho(A,B))$. Suppose that $J$ is a class
of $\G(tB-A)$ such that $\rho(A,B)=\frac{\mu}{1+\mu}$, 
where $\mu=\rho((B_J-A_J)^{-1}A_J)$. Let $m=|J|$. Then $tB-A\in\ls$ with
\[
s \ \left\{ \begin{array}{ll}
      \leq n-1 & \mbox{if $m=n$}\\
      < m & \mbox{if $m<n$.}
                   \end{array}\right.
\] 
\end{theorem}

\pf
That $tB-A\in\ls$ for some $s\in\{0,1,\ldots,n\}$ follows  
from Theorem \ref{zthm}. By Lemma \ref{zlem} (iii), if $t\in(0,\rho(A,B))$, 
then $tB-A\not\in L_n$ since $\rho(A,B)=\tau_n$. Thus $s\leq n-1$. 
When $m<n$, under the assumptions of the theorem, we have
$\tau_n=\rho(A,B)=\frac{\mu}{1+\mu}\leq \tau_{m}$ and hence
$\tau_{m}=\tau_{m+1}=\dots=\tau_n$. It follows that $s<m$.
\eop

We now apply the results of this section to Example \ref{example2},
which has two classes. Class $J=\{2,4\}$ is the only class of 
$\G(A)\cup\G(B)$ such that $(\rho(A,B)B-A)_J$ is singular, and $J$ is 
accessed by no other class. By Theorem \ref{zbasis},
there exists an eigenvector $x$ of $(A,B)$ associated with $\rho(A,B)$
with $x_1=x_3=0$, $x_2>0$ and $x_4>0$. Since $|J|=2$, by Theorem
\ref{zs}, $tB-A\in L_0\cup L_1$ for all $t\in(0,\rho(A,B))$,
agreeing with the exact partition given in Example \ref{example2}.



\begin{thebibliography}{99}

\bibitem{BOV}
R. B. Bapat, D. D. Olesky, and P. van den Driessche.
Perron-Frobenius theory for a generalized eigenproblem.
{\em Linear and Multilinear Algebra}, 40:141-152, 1995.

\bibitem{BP} 
A. Berman and R. J. Plemmons. {\em Nonnegative
Matrices in the Mathematical Sciences},  Academic Press,
New York, 1979.

\bibitem{FM} 
M. Fiedler and T. Markham. A classification of matrices of class Z. 
{\em Linear Algebra and Its Applications}, 173:115-124, 1992.

\bibitem{HJ} 
Roger A. Horn and Charles R. Johnson. {\em Topics in
Matrix Analysis}, Cambridge University Press, 1991.

\bibitem{Na}
Reinhard Nabben. Z-matrices and inverse Z-matrices.
{\em Linear Algebra and Its Applications}, 256:31-48, 1997.


\bibitem{S} 
Hans Schneider.  
The influence of the marked reduced graph of a nonnegative matrix 
on the Jordan Form and on related properties: A survey.
{\em Linear Algebra and Its Applications}, 84:161-189, 1986.

\bibitem{Sm}
Ronald S. Smith. Some results on a partition of Z-matrices.
{\em Linear Algebra and Its Applications}, 223/224:619-629, 1995.

\end{thebibliography}
\end{document}